\documentclass[11pt,a4paper]{article}

\usepackage[utf8]{inputenc}
\usepackage[T1]{fontenc}
\usepackage{lmodern}
\usepackage{geometry}
\usepackage{microtype}
\usepackage{amsmath,amssymb}
\usepackage{graphicx}
\usepackage{booktabs}
\usepackage{tabularx}
\usepackage{array}
\usepackage{caption}
\usepackage{hyperref}
\usepackage{url}
\usepackage{amsthm}

\hypersetup{
  colorlinks=true,
  linkcolor=blue,
  citecolor=blue,
  urlcolor=blue,
  pdftitle={Integer Realization of an Equivelar Octahedron of Genus 3},
  pdfauthor={Ruslan Mizhaev}
}

\graphicspath{{figures/}}

\title{Integer Realization of an Equivelar Octahedron of Genus 3}
\author{Ruslan Mizhaev\\\texttt{\small mimodart@gmail.com}}
\date{September 2026}

\begin{document}
\maketitle

\begin{abstract}
We present an integer-coordinate realization of a genus-3 polyhedral surface with eight planar simple nonagonal faces. It has 24 vertices and 36 edges, and each vertex is incident with three faces. Every pair of faces shares an edge: 20 pairs share one edge and 8 pairs share two. Vertex coordinates, face walks, and plane equations are given. Exact verification confirms that all faces are planar, the surface is closed and orientable, no unintended intersections occur, and the realization has $C_4$ symmetry. The construction is based on an earlier construction and preserves its incidence structure.
\end{abstract}

\section{Introduction}
Equivelar polyhedra provide geometric realizations of regular maps on surfaces of positive genus. They are useful for studying the interaction between combinatorial topology, graph theory, and three-dimensional geometry. Classical examples include the Szilassi polyhedron and related toroidal complexes~\cite{ace,szilassi1986,csaszar1950,grunbaumSzilassi2009}. General existence and classification questions for equivelar maps and polyhedral manifolds have been considered by Datta~\cite{datta2005}, McMullen, Schulz, and Wills~\cite{mcmullenSchulzWills1982}, and Schulte and Wills~\cite{schulteWills2012}. Background on map coloring, regular maps of genus three, and graph embeddings on surfaces is given by Ringel~\cite{ringel1974}, Sherk~\cite{sherk1959}, and Mohar and Thomassen~\cite{moharThomassen2001}.

This paper gives a compact, exactly checkable integer realization of an equivelar octahedron of genus $3$. The underlying surface has type $\{9,3\}$: every face is a nonagon and exactly three faces meet at every vertex. The realization corresponds to variant V2 in the author's earlier work~\cite{mizhaev2020}; the purpose here is to provide a self-contained coordinate certificate that can be reproduced directly in computer algebra or CAD software.

\section{Preliminaries and Definitions}
A polyhedral surface is a finite collection of planar polygonal faces whose intersections are empty, a common vertex, or a common edge, with the local neighborhood of each vertex homeomorphic to a disk. It is \emph{closed} when every edge belongs to exactly two faces, and \emph{orientable} when the faces admit compatible orientations.

An equivelar map of type $\{p,q\}$ is a map in which every face is a $p$-gon and every vertex has valence $q$. Here $p=9$ and $q=3$. The term ``octahedron'' refers to the eight-face combinatorial type; the surface has genus $3$ and is therefore non-convex and non-spherical. We use a face walk $v_1-\cdots-v_9-v_1$ to specify the cyclic order of the vertices of a face.

\section{Main Result}
\label{sec:main}
\textbf{Theorem.} The vertices and face walks in Tables~\ref{tab:vertices} and~\ref{tab:faces}, together with the planar faces defined by them, form a closed orientable equivelar polyhedral surface of type $\{9,3\}$ and genus $3$. The surface has 24 vertices, 36 edges, eight faces, and a geometric symmetry group isomorphic to $C_4$.

\paragraph{Proof.}
The verification is separated into incidence, planarity, topology, and symmetry checks. Each item is an exact calculation over the integers and is summarized in Sections~\ref{sec:planes}--\ref{sec:symmetry}. These checks establish all assertions of the theorem.\qed

\section{Plane Equations and Coordinate Description}
\label{sec:planes}
The eight supporting planes are
\begin{align*}
F_1:&\quad 21x+3y-10z=-1440,\\
F_2:&\quad 21x+3y+10z=1440,\\
F_3:&\quad 3x+z=234,\\
F_4:&\quad 3x-z=-234,\\
F_5:&\quad 3y-z=234,\\
F_6:&\quad 3y+z=-234,\\
F_7:&\quad 3x-21y+10z=-1440,\\
F_8:&\quad 3x-21y-10z=1440.
\end{align*}

\begin{table}[ht]
\centering
\caption{Integer coordinates of the 24 vertices.}
\label{tab:vertices}
\scriptsize
\setlength{\tabcolsep}{5pt}
\begin{tabular}{c|rrr@{\qquad}c|rrr}
\toprule
\multicolumn{4}{c@{\qquad}}{Vertices 1--12} & \multicolumn{4}{c}{Vertices 13--24}\\
\cmidrule(r){1-4}\cmidrule(l){5-8}
No. & $x$ & $y$ & $z$ & No. & $x$ & $y$ & $z$\\
\midrule
1&-72&84&18&13&-9&-147&207\\
2&-36&112&102&14&84&72&-18\\
3&72&-84&18&15&112&36&-102\\
4&36&-112&102&16&-48&-84&18\\
5&0&300&234&17&147&-9&-207\\
6&-84&48&-18&18&-18&126&144\\
7&9&147&207&19&-126&-18&-144\\
8&0&-300&234&20&-300&0&-234\\
9&84&-48&-18&21&-84&-72&-18\\
10&-112&-36&-102&22&18&-126&144\\
11&48&84&18&23&126&18&-144\\
12&-147&9&-207&24&300&0&-234\\
\bottomrule
\end{tabular}
\end{table}

\begin{table}[ht]
\centering
\caption{Cyclic face walks; the last vertex is joined to the first.}
\label{tab:faces}
\begin{tabular}{c l}
\toprule
Face & Vertex walk\\
\midrule
$F_1$ & $6-10-16-22-18-7-5-1-2$\\
$F_2$ & $9-15-11-18-22-13-8-3-4$\\
$F_3$ & $7-5-8-3-17-23-9-15-14$\\
$F_4$ & $13-8-5-1-12-19-6-10-21$\\
$F_5$ & $1-2-11-18-7-14-24-20-12$\\
$F_6$ & $3-4-16-22-13-21-20-24-17$\\
$F_7$ & $14-24-17-23-19-6-2-11-15$\\
$F_8$ & $10-21-20-12-19-23-9-4-16$\\
\bottomrule
\end{tabular}
\end{table}

\section{Geometric and Topological Verification}
\label{sec:verification}
Substitution of every listed vertex into the equation of each incident face gives zero. There are $8\times9=72$ vertex--face incidences, and all 72 equalities hold exactly. For example, vertex 6 satisfies $F_1$ because
\[
21(-84)+3(48)-10(-18)=-1440.
\]

The cyclic walks contain 72 directed edge occurrences. After identifying equal unordered pairs of endpoints, exactly 36 distinct edges remain, and each occurs in precisely two face walks. Every vertex is incident with three edges and three faces; the link of each vertex is a triangle. Thus the surface is a closed 2-manifold.

The Euler characteristic is
\[
\chi=V-E+F=24-36+8=-4,
\]
so orientability gives $\chi=2-2g$ and hence $g=3$. An orientation assignment for $(F_1,\ldots,F_8)$ is
\[
(+,+,-,-,-,-,+,+),
\]
which makes the two induced orientations on every common edge opposite. Hence the surface is orientable.

Pairwise intersection tests between the eight planar polygons show that the intersection of two faces is exactly the edge(s) prescribed by the walks. No additional intersections occur. The multiplicity matrix of common edges is
\[
\begin{pmatrix}
0&1&1&2&2&1&1&1\\
1&0&2&1&1&2&1&1\\
1&2&0&1&1&1&2&1\\
2&1&1&0&1&1&1&2\\
2&1&1&1&0&1&2&1\\
1&2&1&1&1&0&1&2\\
1&1&2&1&2&1&0&1\\
1&1&1&2&1&2&1&0
\end{pmatrix}.
\]
Thus 20 face pairs share one edge and 8 pairs share two edges. In particular, every pair of faces has a common edge. The corresponding cyclic order of the eight faces in the dual multigraph is
\[
F_1-F_4-F_8-F_6-F_2-F_3-F_7-F_5-F_1.
\]

\begin{figure}[ht]
\centering
\includegraphics[width=0.98\textwidth]{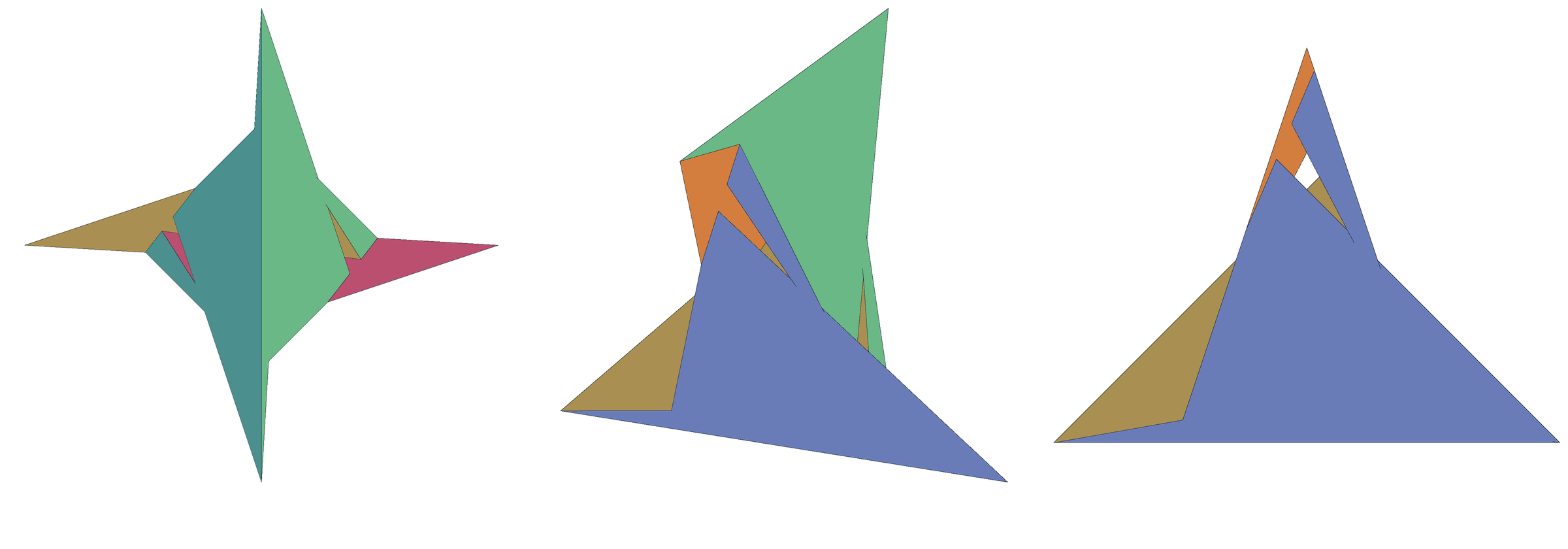}
\caption{Representative views of the integer realization.}
\label{fig:views}
\end{figure}

\begin{figure}[ht]
\centering
\includegraphics[width=0.92\textwidth]{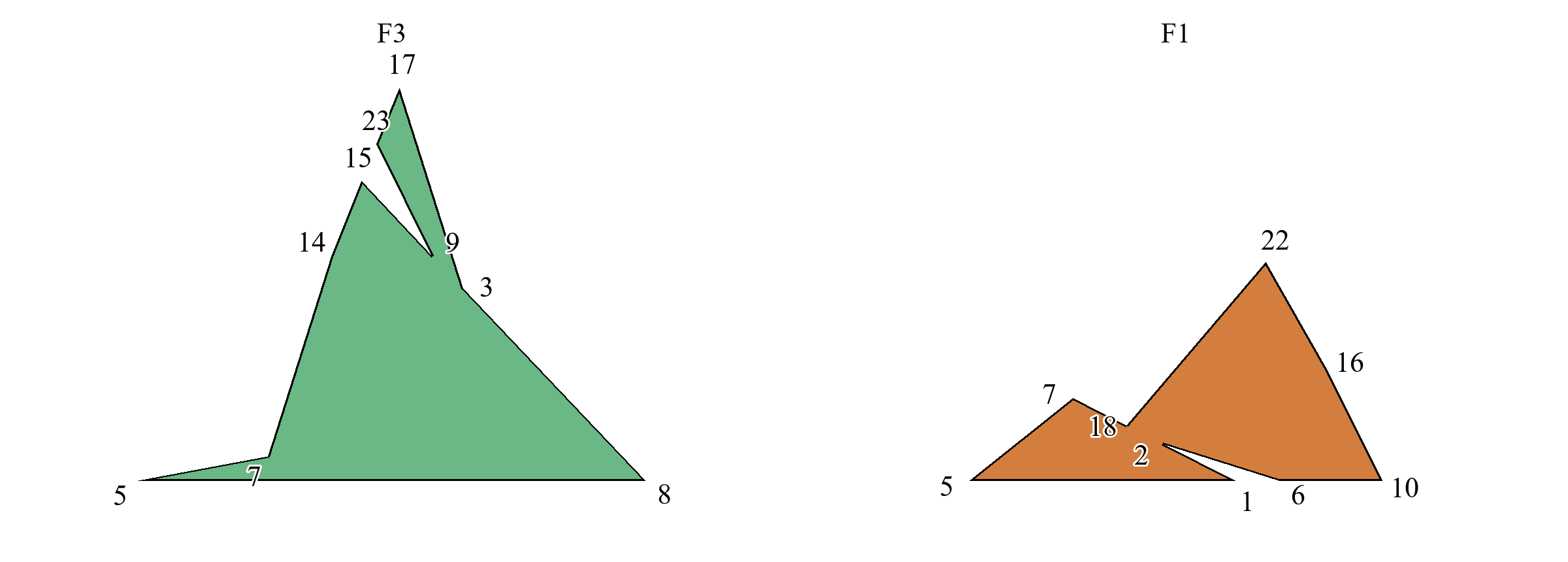}
\caption{Planar representatives of the two face orbits.}
\label{fig:faces}
\end{figure}

\section{Symmetry and Adjacency}
\label{sec:symmetry}
The linear map
\[
T(x,y,z)=(y,-x,-z)
\]
preserves the vertex set, the face planes, and the face walks. Since $T^4=\mathrm{id}$, it generates a cyclic symmetry group $\langle T\rangle\cong C_4$. Its action on the faces is
\[
(F_1\ F_7\ F_2\ F_8)(F_3\ F_6\ F_4\ F_5),
\]
so the faces form two orbits of four. Exhaustive enumeration of the automorphisms of the incidence graph gives four automorphisms, and every one is induced by a power of $T$; hence the geometric symmetry group is exactly $C_4$.

The dual multigraph has eight vertices and 36 edges. Every dual vertex has degree 9 when common edges are counted with multiplicity, while its underlying simple graph is $K_8$. This records the fact that every pair of faces is adjacent, with eight pairs adjacent twice.

\section{Reproducibility and Limitations}
All data needed to reproduce the realization are contained in Tables~\ref{tab:vertices} and~\ref{tab:faces} and in the eight plane equations above. The checks can be implemented using exact integer arithmetic; no floating-point tolerance is required. The certificate verifies the stated realization and its incidence structure. It does not establish coordinate minimality, optimal aspect ratios, or historical priority.

\section*{Acknowledgments}
The author acknowledges the use of OpenAI's ChatGPT for text translation, academic formatting, and assistance with certain calculations. The author takes full intellectual responsibility for the correctness and content of this work.

\bibliographystyle{unsrt}
\bibliography{references}

\end{document}